\documentstyle [twoside,12pt]{article}
\font\twelvebb=msym10 at 12pt
\font\twelvegoth=eufm10 at 12pt
\newfam\gothfam
\textfont\gothfam=\twelvegoth

\newfam\bbfam
\textfont\bbfam=\twelvebb
\def\bb{\fam\bbfam\twelvebb}

\newcommand{\nc}{\newcommand}
\nc{\si}{\sigma}\nc{\al}{\alpha}\nc{\be}{\beta}
\nc{\la}{\lambda}
\nc{\th}{\theta}    \nc{\dst}{\displaystyle}       
\nc{\ga}{\gamma}     \nc{\nnb}{\nonumber}           \nc{\nf}{\infty}
\nc{\Ga}{\Gamma}
\nc{\de}{\delta}               \nc{\sn}{{\rm sn\,}}   \nc{\sd}{{\rm sd\,}}
\nc{\eps}{\epsilon}            \nc{\cn}{{\rm cn\,}}   \nc{\nd}{{\rm nd\,}}
\nc{\ra}{\rightarrow}          \nc{\dn}{{\rm dn\,}}
\nc{\Lra}{\Longrightarrow}                             \nc{\ssc}{{\rm sc\,}}
\nc{\beq}{\begin{equation}}                            \nc{\cd}{{\rm cd\,}}
\nc{\eeq}{\end{equation}}                              \nc{\dc}{{\rm dc\,}}
\nc{\beqa}{\begin{eqnarray}}
\nc{\eeqa}{\end{eqnarray}}
\nc{\prh}[2]{\left(\begin{array}{c}#1 \\#2 \end{array}\right)}
\nc{\prf}[3]{\left(\begin{array}{c}#1 \\#2 \end{array}\right.\left.
;\rule{0 cm}{0.6cm}#3\right)}
\nc{\cro}[2]{\left[\begin{array}{c}#1 \\#2
\end{array}\right]}
\nc{\fx}[5]{{_{#1}F_{#2}}\prf{#3}{#4}{#5}}
\nc{\COS}{\int_0^{\th(w)}du\ \cos2\sqrt{\si}(\th(w)-u)}
\nc{\SIN}{\int_0^{\th(w)}du\ \frac{\sin2\sqrt{\si}(\th(w)-u)}{2\sqrt{\si}}}
\nc{\SC}{Stieltjes-Carlitz }

\title{ \bf Associated Stieltjes-Carlitz polynomials and a generalization
of Heun's differential equation.}
\author{Galliano VALENT\thanks
{\noindent Laboratoire de Physique Th\'eorique et Hautes Energies,
 Unit\'e associ\'ee au CNRS UA 280,~Universit\'e PARIS 7,
 2 Place Jussieu, F-75251 CEDEX 05.}}
\date{}
\begin{document}
\maketitle
\vspace{-8cm}\centerline{ COMMUNICATION AT THE 4th INTERNATIONAL CONFERENCE}
\centerline{ ON ORTHOGONAL POLYNOMIALS. EVIAN. OCTOBER 1992.}\vspace{+8cm}
\begin{abstract}
\noindent The generating function of Stieltjes-Carlitz polynomials is a
solution of Heun's differential equation and using this relation Carlitz was
the first
to get exact closed forms for some Heun functions. Similarly the associated
Stieltjes-Carlitz polynomials lead to a new differential equation which we call
associated Heun.
Thanks to the link with orthogonal polynomials we are able to deduce two
integral relations connecting associated Heun functions with different
parameters and to
exhibit the set of associated Heun functions which generalize Carlitz's. Part
of these
results were used by the author to derive the Stieltjes transform of the
measure of
orthogonality for the associated Stieltjes-Carlitz polynomials using asymptotic
analysis; here we present a new derivation of this result.
\end{abstract}

\section{Introduction}
The Stieltjes-Carlitz polynomials are othogonal polynomials with a three term
recurrence
relation
\beq\label{a1}\begin{array}{l}
(\lambda_n+\mu_n-x)F_n(x)=\mu_{n+1}F_{n+1}+\lambda_{n-1}F_{n-1},
\hspace{1cm} n\geq 0,\\ [0.2cm]
F_{-1}(x)=0,\hspace{2cm} F_0(x)=1,\end{array}\eeq
with the rates
\beq\label{a2}\left\{\begin{array}{l}
\la_n=k^2(2n+1)^2,\hspace{1.5cm}\mu_n=(2n)^2\\ [0.2cm]
\la_n=(2n+1)^2,\hspace{1.5cm}\mu_n=k^2(2n)^2\end{array}\right.\eeq
and $0<k^2<1$. Their orthogonality measure, given for instance in \cite[p.
194]{ch78},
 was first derived by Stieltjes \cite{st94}
 using continued fraction techniques. Later on Carlitz  \cite{ca60} obtained a
generating function for the $F_n(x)$ from
which he was able to derive anew the orthogonality measure. The most striking
fact, which
does not seem to have been realized by Carlitz himself, is that the generating
functions he
had obtained give a finite set of exact solutions of Heun's differential
equation. The
complete list of Carlitz results were gathered in  \cite{val86}, and an
integral
transformation connecting Heun functions with different parameters was derived.
This last
tool, combined with Carlitz results, led to an enlarged set of Heun functions.

More recently we have obtained a generating function for the associated
Stieltjes-Carlitz
polynomials with rates
\beqa\label{a3}
\la_n&&\hspace{-0.6cm}=k^2(2n+2c+1)^2,
\hspace{2cm}\mu_n=4(n+c)^2+\mu\de_{n0}\\[0.2cm]
\label{a4}
\la_n&&\hspace{-0.6cm}=(2n+2c+1)^2,
\hspace{2cm}\mu_n=4k^2(n+c)^2+k^2\mu\de_{n0}\eeqa
This result, combined with asymptotic analysis, has led to the Stieltjes
transform of the
associated Stieltjes-Carlitz polynomials \cite{val94}

Our aim is to show that these results lead to a finite set of exact solutions
of what
could be called the ``associated Heun" differential equation which emerges as
an equation
satisfied by the generating functions of the polynomials $F_n(x)$ with the
rates
(\ref{a3}), (\ref{a4}).

This link with orthogonal polynomials is even more fruitful since it gives a
convenient
tool to derive two new integral transformations relating associated Heun
functions with
different parameters.

The plan of this communication is the following.

Section 2 is devoted to a short summary on Heun's differential equation and its
relation
with orthogonal polynomials.

In section 3 we present the associated Heun differential equation and derive
two
integral connection relations for its solutions.

In section 4 we give a finite set of exact solutions of the associated Heun
differential
equation which generalize  Carlitz results to non-vanishing $(c,\mu)$.

In section 5, in order to cross-check the Stieltjes transform obtained in
\cite{val94},
 we present a completely different derivation which uses Karlin and MacGregor
representation theorem and the link between orthogonal polynomials and birth
and death
processes. Solving the Kolmogorov equations leads to a Stieltjes transform
which is
in perfect agreement with asymptotic analysis.

\section{Heun differential equation}
\noindent This equation is the most general second order differential equation
with four
regular singular points
$$\hspace{3cm}w=0,\ 1,\ 1/k^2,\ +\nf\hspace{2cm}0\leq k^2\leq 1$$
and is described by the array of parameters
$$P=\{\al,\be;\ga,\de,\eps;s\}\hspace{2cm}\al+\be=\ga+\de+\eps-1.$$
The accessory parameter $s$ is unconstrained. We shall follow the standard
notations of
\cite[p.576]{ww65},\cite[p.57-62]{emot3} for Heun's differential equation
\beqa\nnb
w(1-w)(1-k^2w)D^2 F+[\ga(1-w)(1-k^2w)-\de w&&\hspace{-0.6cm}(1-k^2w)-\eps
k^2w(1-w)]DF
\\ [0.2cm]\label{b1}
+&&\hspace{-0.6cm}(\al\be k^2w+s)F=0\eeqa
with $\dst D=\frac{d}{dw}$.

We shall denote by $Hn(P,w)$ the unique solution of (\ref{b1}) which is
analytic for
$|w|<1$ and is normalized according to
\beq\label{bb1}
Hn(P,w=0)=1.\eeq

This equation, for arbitrary values of $s$, can be solved in terms of
hypergeometric
functions only for two values of the parameter $k^2$:

1) if $k^2=0$ the solution analytic around $w=0$ is
$$\dst\fx{2}{1}{r_{+},r_{-}}{\ga}{w}\hspace{1cm}r_{\pm}=a\pm \sqrt{a^2+s}
\hspace{1cm}a=\frac{\ga+\de-1}{2}$$

2) if $k^2=1$ the solution analytic around $w=0$ was shown in \cite{rv82} to be
$$(1-w)^r\fx{2}{1}{r+\al,r+\be}{\ga}{w}\hspace{1cm}r=a+\sqrt{a^2-\al\be
-s},\hspace{1cm}
a=\frac{\ga -\al-\be}{2}$$

For some particular values of the accessory parameter $s$ Heun's functions
degenerate into
hypergeometric functions of the variable $R(w)$, where $R(w)$ is a polynomial
of second
degree in $w$. These values of $s$ are listed in \cite{ku79}.

In all what follows we shall not consider these particular cases.

Since $Hn(P,w)$ is analytic for $|w|<1$ we can consider it as the generating
function of
the polynomials $F_n(P,s)$, with variable $s$, such that
\beq\label{b2}
Hn(P,w)=\sum_{n\geq 0}F_n(P,s)w^n\hspace{2cm}|w|<1\eeq
Relations (\ref{b1},\ref{b2}) imply routinely the three term recurrence
relation for the $F_n$
\beq\label{b3}\left.\begin{array}{c}
 (\la_n+\mu_n+\ga_n-s-\al\be k^2)F_n=\mu_{n+1}F_{n+1}+\la_{n-1}F_{n-1},
\hspace{1cm}n\geq 0\\ [0.2cm]F_{-1}=0,\hspace{1cm}F_0=1\end{array}\right\}\eeq
with
$$\la_n=k^2(n+\al)(n+\be),
\hspace{1cm}\mu_n=n(n+\ga-1),\hspace{1cm}\ga_n=(1-k^2)\de n.$$
This recurrence exhibits the polynomial character of $F_n$ with respect to
either the
variable $s$ or the more familiar $x=s+\al\be k^2$.

One should observe on (\ref{b3}) that only for $\de=0$ do we have a true birth
and death
process with birth rate $\la_n$ and death rate $\mu_n$ (see \cite{ilmv90} for
an
introduction). For $\de\neq 0$ we have {\em killing} in the sense of Karlin and
Tavar\'e
\cite{kt82},\cite{kt82a} with rate $\ga_n$.

The finite set of exact solutions of (\ref{b1}) obtained by Carlitz can be
found in \cite[p.692]{val86} as well as as the following integral transform,
quoted here
for convenience.

Let us define
$P'=\{\al'=\ga,\be'=\be;\ga'=\al,\de'=\de+\ga-\al,\eps'=\eps+\ga-\al;s\}$,
then provided that ${\rm Re}\,\ga>{\rm Re}\,\al>0$, we have
\beq\label{b4}
Hn(P,w)=\dst\frac 1{B(\al,\ga-\al)}\int_0^1dt\ t^{\al-1}(1-t)^{\ga-\al-1}
Hn(P';wt)\eeq
for  $w$ in ${\bb C}\backslash [1,+\nf [.$ This relation, combined with the set
of
Carlitz solutions gives another finite set of solutions described in
\cite[p.693]{val86}.

\section{Integral connection relations for associated Heun functions}
\noindent Let us turn ourselves to the associated polynomials with the
recurrence
\beq\label{c1}\begin{array}{c}
\left(\la_n+\mu_n+\ga_n-s-(\al+c)(\be+c) k^2+k^2\de
c\right)F_n=\mu_{n+1}F_{n+1}+\la_{n-1}F_{n-1},
\hspace{0.2cm}n\geq 0\\ [0.2cm]F_{-1}=0,\hspace{1cm}F_0=1\end{array}\eeq
with
$$\left\{\begin{array}{l}
\la_n=k^2(n+c+\al)(n+c+\be)\\ [0.2cm]
\mu_n=(n+c)(n+c+\ga-1)+\mu\de_{n0}\\ [0.2cm]
\ga_n=(1-k^2)\de(n+c)\end{array}\right.$$
Two new parameters appear: $c$ which is an association parameter, and $\mu$
which is a
co-recursivity parameter. Clearly, for vanishing $(c,\mu)$ the recurrence
(\ref{c1})
reduces to (\ref{b3}). The constant terms added to $s$ were chosen for
notational
convenience reasons.

We define
\beq\label{c2}
Hn(c,\mu,P;w)=\sum_{n\geq 0}F_n(c,\mu,P)w^n\hspace{2cm}|w|<1\eeq
and from (\ref{c1}) it is easy to obtain
\beqa\nnb
w(1-w)&&\hspace{-0.6cm}(1-k^2w)D^2F+[(\ga+2c)(1-w)(1-k^2w)-\de w(1-k^2w)-\eps
k^2w(1-w)]DF\\ [0.2cm]\label{c3}\dst
+&&\hspace{-0.6cm}\left[(\al+c)(\be+c)k^2w+\frac{c(c+\ga-1)}{w}(1-w)+s-\de
c\right]F=
\frac{c(c+\ga-1)}{w}+\mu\eeqa
In view of its origin it is natural to call this equation the associated Heun
differential
equation. This departs from the usual terminology where associated differential
equations
refer to {\em homogeneous} extensions of a given differential equation whereas
here we
have an {\em inhomogeneous} extension.

At any rate it reduces to Heun equation in the particular cases $(c=0,\mu=0)$
and
 $(c=1-\ga,\mu=0)$, according to the relations
\beq\label{cc3}\begin{array}{l}
Hn(0,0,P;w)=Hn(P;w)\\ [0.2cm]
Hn(1-\ga,0,P;w)=Hn(\tilde{P};w)\end{array}\eeq
with the array $\tilde{P}=\{1-\ga+\al,1-\ga+\be;2-\ga,\de,\eps;s-(1-\ga)\de\}.$

Let us now derive two integral transformations relating associated Heun
functions. To do
this we switch from the $F_n$ to the $G_n$ defined by
\beq\label{c4}
G_0=F_0,\hspace{2cm}G_n=\mu_1\cdots\mu_nF_n=
(1+c)_n(\ga+c)_nF_n\hspace{1cm}n\geq 1\eeq
whose recurrence
\beq\label{c5}\begin{array}{c}
\left(\la_n+\mu_n+\ga_n-s-(\al+c)(\be+c) k^2+k^2\de
c\right)G_n=G_{n+1}+\la_{n-1}\mu_{n}G_{n-1},
\hspace{1cm}n\geq 0\\ [0.2cm]G_{-1}=0,\hspace{1cm}G_0=1\end{array}\eeq
reveals that $(-1)^n G_n$ is monic in the variable
$x=s+k^2(\al+c)(\be+c)-k^2\de c$.

The basic technique to get an integral transform is to look for a mapping of
the
parameters P which leaves invariant the recurrence (\ref{c5}).

A first possibility is the mapping $P'_{\al}$
$$\begin{array}{l}\al'=\ga,\ \be'=\be\\ [0.2cm]
\ga'=\al,\ \de'=\de+\ga-\al,\ \eps'=\eps+\ga-\al\\ [0.2cm]
\ s'=s,\ c'=c,\ \mu'=\mu\end{array}$$
for which we have
$$G_n(P'_{\al})=G_n(P)\hspace{2cm}n\geq 0.$$
Using (\ref{c4}) gives
\beq\label{c6}
\dst F_n(P)=\frac{(c+\al)_n}{(c+\ga)_n}F_n(P'_{\al}).\eeq
If ${\rm Re}\ \ga>{\rm Re}\ \al>-{\rm Re}\ c$ we can write
$$\dst\frac{(c+\al)_n}{(c+\ga)_n}=\frac 1{B(\ga-\al,\al+c)}\int_0^1dt\
t^{n+c+\al-1}
(1-t)^{\ga-\al-1}$$
and inserting this in (\ref{c6}), multiplying each term by $w^n$ and summing
$n$ from zero
to infinity gives
\beq\label{c7}
Hn(c,\mu,P,w)=\frac 1{B(\ga-\al,\al+c)}\int_0^1dt\ t^{c+\al-1}
(1-t)^{\ga-\al-1}Hn(c,\mu,P'_{\al};wt)\eeq
valid for $|w|<1$. The term by term integration is allowed since the right hand
side power
series is absolutely and uniformly convergent for $|w|\leq R<1$. Analytic
continuation extends this relation to ${\bb C}\backslash [1,+\nf[.$ Clearly for
$c=\mu=0$ we
recover (\ref{b4}) and there is another integral transformation $P'_{\be}$
obtained from
$P'_{\al}$ by the exchange of the couples $(\al,\al')$ and $(\be,\be')$.

A second possibility which leaves invariant the recurrence (\ref{c5}) is
$P''_{\al}$ with
$$\begin{array}{l}\al''=2-\al,\ \be''=\be+1-\al\\ [0.2cm]
\ga''=\ga+1-\al,\ \de''=\de+1-\al,\ \eps''=\eps+1-\al\\ [0.2cm]
s''=s+(\al-1)(\ga+\de-\al),\ c''=c+\al-1,\ \mu''=\mu\end{array}$$
which leads to
$$\dst F_n(P)=\frac{(c+\al)_n}{(c+1)_n}F_n(P''_{\al}).$$
The corresponding integral transform follows analogously to (\ref{c7})
$$Hn(c,\mu,P,w)=\frac 1{B(1-\al,c+\al)}\int_0^1dt\ t^{c+\al-1}
(1-t)^{-\al}Hn(c+\al-1,\mu,P''_{\al};wt)$$
and is valid for $1>{\rm Re}\ \al>-{\rm Re}\ c$ and $ w\in {\bb C}\backslash
[1,+\nf[.$

This second integral relation is a genuinely new result, since it changes the
value of the
association parameter from $c$ to $c+\al-1$. For this reason it could not
appear in the
previous analyses where $c=0$.
Here too the interchange of the couples  $(\al,\al'')$ and $(\be,\be'')$ leads
to another
mapping $P''_{\be}$.

\section{Exact solutions of the associated Heun equation}
\noindent Before giving the set of associated Heun functions which generalize
Carlitz ones
we shall explain, on the first of them, how they can be constructed.

We first make the change of function
\beq\label{d1}G=w^cF\eeq
which brings (\ref{c3}) to
$$\begin{array}{l}w(1-w)(1-k^2w)D^2 G+[\ga(1-w)(1-k^2w)-\de w(1-k^2w)-\eps
k^2w(1-w)]DG
\\ [0.2cm]
\hspace{3cm}+\left[s+k^2c(c+\eps+\ga-1)+\al\be
k^2w\right]G=c(c+\ga-1)w^{c-1}+\mu w^c.
\end{array}$$
The first exact solution will correspond to the parameters
$$P=\left\{\al=0,\ \be=\frac 12;\ \ga=\frac 12,\ \de=\frac 12,\ \eps=\frac 12;
\ s=\si-k^2c^2\right\}\hspace{1cm}c>0.$$
The change of variable
\beq\label{d2}
\sqrt{w}=\sn(\th;\ k^2)\eeq
(in what follows, concerning elliptic functions we stick to the notations of
\cite{ww65};
here $\sqrt{w}$ is the square root which is positive for real positive $w$)
reduces the
differential equation to
$$\partial^2_{\th}G+4\si G=2c(2c-1)(\sn^2\th)^{c-1}+4\mu(\sn^2\th)^c=J(\th)$$
Its solution is
\beq\label{d3}\dst
G(\th)=\int_0^{\th}du\ \frac{\sin2\sqrt{\si}(\th-u)}{2\sqrt{\si}}J(u)\eeq
and is meaningful provided that $c>1/2$. In order to extend this integral
representation
to $c>0$ an integration by parts of $(\sn^2u)^{c-1}$ is needed with the final
result
\beqa\nnb\dst
w^c&&\hspace{-0.6cm}Hn(c,\mu,P;\si-k^2c^2;w)=\COS 2c(\sn^2u)^{c-1/2}\cn u\dn
u\\ \label{d4}
+&&\hspace{-0.6cm}\SIN
\left[4(c^2+c^2k^2+\mu)(\sn^2u)^c-2c(2c+1)k^2(\sn^2u)^{c+1}\right]
\eeqa
The variable $\th(w)$ is obtained through the inversion of relation (\ref{d2})
$$\dst\th(w)=\int_0^{\sqrt{w}}\frac{dt}{\sqrt{(1-t^2)(1-k^2t^2)}}.$$
It is analytic for $|w|<1$; its analytic extension to the whole complex $w$
plane has
been described in \cite[p.122]{ak2} and is analytic but for the branch
points $w=1,\ 1/k^2.$

We take for $(\sn^2u)^{c-1/2}$ and $w^c$ one and the same principal branch in
order to
secure the analyticity of $Hn(w)$ for $|w|<1$.

It is lengthy, even if straightforward, to check that (\ref{d4}) is indeed a
solution of
(\ref{c3}) in the complex plane deprived with the points $w=1,\ 1/k^2$ and that
the
normalization condition (\ref{bb1}) does hold.

As mentioned in the previous section, this associated Heun function should
reduce to a
Heun function either if $(c=\mu=0)$ or if $(c=1/2,\mu=0)$. In the first case a
limiting
procedure which makes use of
$$\lim_{c\ra 0}\int_0^{\th}du\ f(u)\ 2c(\sn^2u)^{c-1/2}=f(0)$$
gives
$$\dst\lim_{c\to 0}Hn(c,\mu,P;\si;w)=\cos2\sqrt{\si}\th(w) +\frac{\mu}{\si}
\left(1-\cos2\sqrt{\si}\th(w)\right).$$
In the second limiting case, using relation (\ref{d3}) we get
$$\dst\lim_{c\ra 1/2}\sqrt{w}Hn\left(c,\mu,P;\si-\frac{k^2}{4};w\right)=
\frac{\sin2\sqrt{\si}\th(w)}{2\sqrt{\si}}+4\mu\SIN\ \sn u$$
For $\mu=0$, using the relations (\ref{cc3}) we recover Carlitz results:
$$\left\{\begin{array}{l}\dst
Hn\left(0,\frac 12;\frac 12,\frac 12,\frac
12;s=\si;w\right)=\cos(2\sqrt{\si}\th(w))
\\ [0.4cm]\dst
\sqrt{w}Hn\left(\frac 12,1;\frac 32,\frac 12,\frac 12;s=\si-\frac
{1+k^2}4;w\right)=
\frac{\sin2\sqrt{\si}\th(w)}{2\sqrt{\si}}\end{array}\right.$$
In order to get the remaining set of exact solutions one has to change
(\ref{d1}) into
$$G=w^{c+\la}(1-w)^{\mu}(1-k^2w)^{\nu}F$$
where $\la,\mu,\nu$ take the values $0$ or $1/2$. We get in this way seven more
solutions
to be listed below.

$\bullet$
$\dst P=\left\{\frac 12,1;\frac 12,\frac 32,\frac 12;s=\si-\frac
14-k^2c^2\right\}$

\beq\label{d5}\dst\begin{array}{l}\dst
w^c\sqrt{1-w}Hn(c,\mu,P;w)=\COS\ 2c(\sn^2u)^{c-1/2}\dn u\\ [0.2cm]\dst
\hspace{5cm}+4(k^2c^2+\mu)\SIN\ (\sn^2u)^c\cn u\end{array}\eeq
which reduces for $(c=\mu=0)$ and $(c=1/2,\mu=0)$ to
$$\left\{\begin{array}{l}\dst
\ \sqrt{1-w}Hn\left(\frac 12,1;\frac 12,\frac 32,\frac 12;s=\si-\frac 14
;w\right)=
\cos(2\sqrt{\si}\th(w))\\ [0.4cm]\dst
\sqrt{w(1-w)}Hn\left(1,\frac 32;\frac 32,\frac 32,\frac
12;s=\si-1-\frac{k^2}{4};w\right)=
\frac{\sin2\sqrt{\si}\th(w)}{2\sqrt{\si}}\end{array}\right.$$

$\bullet$
$\dst P=\left\{\frac 12,1;\frac 12,\frac 12,\frac 32;s=\si-k^2(c+\frac
12)^2\right\}$

\beq\label{d6}\dst\begin{array}{l}\dst
w^c\sqrt{1-k^2w}Hn(c,\mu,P;w)=\COS\ 2c(\sn^2u)^{c}\cn u\\ [0.2cm] \dst
\hspace{5cm}+4(c^2+\mu)\SIN\ (\sn^2u)^c\dn u\end{array}\eeq
which reduces for $(c=\mu=0)$ and $(c=1/2,\mu=0)$ to
$$\left\{\begin{array}{l}\dst
\ \sqrt{1-k^2w}Hn\left(\frac 12,1;\frac 12,\frac 12,\frac
32;s=\si-\frac{k^2}{4};w\right)=
\cos(2\sqrt{\si}\th(w))\\ [0.4cm]\dst
\sqrt{w(1-k^2w)}Hn\left(1,\frac 32;\frac 32,\frac 12,\frac 32;s=\si-\frac
14-k^2;w\right)=
\frac{\sin2\sqrt{\si}\th(w)}{2\sqrt{\si}}\end{array}\right.$$

$\bullet$
$\dst P=\left\{1,\frac 32;\frac 12,\frac 32,\frac 32;s=
\si-\frac 14-k^2(c+\frac 12)^2\right\}$

$$\dst\begin{array}{l}\dst
\sqrt{(1-w)(1-k^2w)}Hn(c,\mu,P;w)=\COS\ 2c(\sn^2u)^{c-1/2}\\ [0.2cm]\dst
\hspace{5cm}+4\mu\SIN\ (\sn^2u)^c\cn u\dn u\end{array}$$
which reduces for $(c=\mu=0)$ and $(c=1/2,\mu=0)$ to
$$\left\{\begin{array}{l}\dst
\ \sqrt{(1-w)(1-k^2w)}Hn\left(1,\frac 32;\frac 12,\frac 32,\frac 32;s=
\si-\frac{1+k^2}4;w\right)=\cos(2\sqrt{\si}\th(w))\\ [0.4cm]\dst
\sqrt{w(1-w)(1-k^2w)}Hn\left(\frac 32,2;\frac 32,\frac 32,\frac
32;s=\si-1-k^2;w\right)=
\frac{\sin2\sqrt{\si}\th(w)}{2\sqrt{\si}}\end{array}\right.$$

The eight particular cases for which either $(c=\mu=0)$ or $(c=1/2,\mu=0)$
reproduce the
results collected in the table \cite[p.692]{val86}. There are four other
solutions:

$\bullet$
$\dst P=\left\{\frac 12,1;\frac 32,\frac 12,\frac 12;s=
\si-\frac 14-k^2(c+\frac 12)^2\right\}$
\newline Let us define
$$F(u)=2c(2c+1)(\sn^2u)^{c-1/2}+4\mu(\sn^2u)^{c+1/2}$$ then we have
$$\sqrt{w}Hn(c,\mu;w)=\SIN\ F(u)$$

$\bullet$
$\dst P=\left\{1,\frac 32;\frac 32,\frac 32,\frac 12;s=
\si-1-k^2(c+\frac 12)^2\right\}$
$$\sqrt{w(1-w)}Hn(c,\mu;w)=\SIN\ F(u)\cn u$$

$\bullet$
$\dst P=\left\{1,\frac 32;\frac 32,\frac 12,\frac 32;s=\si-\frac
14-k^2(c+1)^2\right\}$
$$\sqrt{w(1-k^2w)}Hn(c,\mu;w)=\SIN\ F(u)\dn u$$

$\bullet$
$\dst P=\left\{\frac 32,2;\frac 32,\frac 32,\frac
32;s=\si-1-k^2(c+1)^2\right\}$
$$\sqrt{w(1-w)(1-k^2w)}Hn(c,\mu;w)=\SIN\ F(u)\cn u\dn u$$

For these last cases the limiting cases $(c=0,\mu=0)$ and $(c=1/2,\mu=0)$ do
not give
anything new. Furthermore the reader can check that all the solutions given
here are
correctly normalized at $w=0$.

The results (\ref{d5},\ref{d6}) were first derived in \cite{val94} whilst all
the other are
new. Obviously this set of solutions can be further enlarged using the two
integral
transforms of section 3.

\section{Associated Stieltjes-Carlitz polynomials}
\noindent In \cite{val94} the Stieltjes transform of the orthogonality measure
of the
associated Stieltjes-Carlitz polynomials has been derived for the the first
time. The main
tool used in this work is Markov theorem and asymptotic analysis. It is of some
interest
to check this result using a completely different approach and this is the aim
of this
section.

The strategy used here was already applied to the \SC polynomials in
\cite{val92} and
provided for a new derivation of the orthogonality measures. Its
generalization, described in \cite{val91} and in \cite{val92}, led to a one
parameter family of orthogonality
measures for
the indeterminate moment problem corresponding to the rates
$$\la_n=(4n+1)(4n+2)^2(4n+3),\hspace{1cm}\mu_n=(4n-1)(4n)^2(4n+1).$$

What is basic in this approach is the connection between orthogonal polynomials
and birth
and death processes; the whole problem to get the orthogonality measure is
reduced to the
resolution of a linear partial differential equation.

Let us first recall that the first family of associated \SC polynomials are
defined by the
recurrence relation
\beq\label{e1}\begin{array}{c}
 (\la_n+\mu_n-x)F_n=\mu_{n+1}F_{n+1}+\la_{n-1}F_{n-1},
\hspace{1cm}n\geq 0\\ [0.2cm]F_{-1}=0,\hspace{1cm}F_0=1\end{array}\eeq
with the rates
$$\la_n=k^2(2n+2c+1)^2,\hspace{1cm}\mu_n=4(n+c)^2+\mu\de_{n0}
\hspace{1cm}c\geq 0$$

In order to obtain the Stieltjes transform of their orthogonality measure we
shall relate
them to the birth and death process whose Kolmogorov equation is
\beq\label{e3}\begin{array}{c}
\dst\frac{d}{dt}{\cal P}_{m,n}(t)=\lambda _{n-1}{\cal P}_{m,n-1}(t)
+\mu _{n+1}{\cal P}_{m,n+1}(t)
-(\lambda _{n}+\mu _n){\cal P}_{m,n}(t)\\ [0.2cm]
{\cal P}_{m,n}(0)=\delta_{m n}.\end{array}\eeq
${\cal P}_{m,n}(t)$ is the probability of a population $n$ at time $t$ provided
that it
was $m$ at time $t=0$. It is therefore positive and bounded by one.
The link between (\ref{e1}) and (\ref{e3}) is provided by Karlin and McGregor
representation theorem \cite{km55},\cite{km57}
$${\cal P}_{m,n}(t)=\frac{1}{\pi_m}\int_0^{\infty}d\Psi(x)F_m(x)F_n(x)e^{-tx}$$
with
$$\pi_0=1,\hspace{0.5cm}\pi_m=\frac{\la_0\ldots\la_{m-1}}{\mu_1\ldots\mu_m},
\ \ m=1,2,\ldots$$
{}From this representation theorem it follows that a possible way to get $\Psi$
is to
compute the Laplace transform of
$${\cal P}_{00}(t)=\int_0^{\infty}d\Psi(x)\ e^{-tx}$$
which we shall write
$$\tilde{{\cal P}}_{00}(p)=\int_0^{\infty}dt\ e^{-pt}{\cal P}_{00}(t).$$
For ${\rm Re}\ p>0$ this is nothing but
$$\dst\tilde{{\cal P}}_{00}(p)=\int_0^{\infty}\frac{d\Psi(x)}{p+x}$$
closely related to the Stieltjes transform of the orthogonality measure since
we have
$$\tilde{{\cal P}}_{00}(p)=-S(-p).$$
In this approach no recourse to asymptotic analysis is needed to get $S(z)$: we
just
require $\tilde{{\cal P}}_{00}(p)$. We shall describe in the following how this
can be
worked out just by solving linear partial differential equations.

As a first step we consider the change of basis ${\cal P}_{m,n}(t)\ra
P_{mn}(t)$ such that
\beq\label{ee4}\dst
P_{mn}(t)=\frac{(1+c)_n}{(1/2+c)_n}{\cal P}_{m,n}(t)\eeq
Kolmogorov equation becomes
\beq\label{e4}\begin{array}{l}\dst
\frac{d}{dt}P_{m,n}(t)=\tilde{\la}_{n-1}P_{m,n-1}(t)
+\tilde{\mu}_{n+1}P_{m,n+1}(t)-(\la_{n}+\mu _n)P_{m,n}(t)\\ [0.2cm]\dst
P_{m,n}(0)=\frac{(1+c)_m}{(1/2+c)_m}\delta_{m n}.\end{array}\eeq
with
$$\tilde{\la}_n=k^2(2n+2c+1)(2n+2c+2),
\hspace{1cm}\tilde{\mu}_n=(2n+2c-1)(2n+2c)$$
whilst $\la_n,\mu_n$ are given by (\ref{a3}).

In order to solve (\ref{e4}) we introduce the generating function
$$H_m(t,w)=\sqrt{1-k^2w}\sum_{n\geq 0}P_{mn}(t)\ w^{n+c}\hspace{1cm}|w|<1.$$
The Kolmogorov equation becomes a linear partial differential equation for
$H_m$
\beq\label{e5}\begin{array}{l}
\partial_t
H_m(t,w)=\left\{4w(1-w)(1-k^2w)\partial_w^2+2[(1-w)(1-k^2w)-w(1-k^2w)\right.
\\ [0.2cm]
\hspace{3cm}\left.-k^2w(1-w)]\partial_w\right\}H_m(t,w)
-\left(\mu+\frac{2c(2c-1)}{w}\right)w^c\sqrt{1-k^2w}P_{m0}(t)\end{array}\eeq
with the boundary condition
$$H_m(0,w)=\dst\frac{(1+c)_m}{(1/2+c)_m}w^{m+c}\sqrt{1-k^2w}.$$
It is convenient, from a notational point of view, to keep $P_{m0}(t)$; however
this is
related to the generating function $H_m(t,w)$ by
$$P_{m0}(t)=\lim_{w\to 0}w^{-c}H_m(t,w)$$
{}From now on we shall restrict ourselves to $m=0$, and in order to simplify
(\ref{e5}) we
change the variable to $w=\sn^2(\th,k^2)$ which maps $[0,1]$ into $[0,K]$.
Deleting the
subscript $m=0$ in $H_0$ we are led to
\beq\label{e6}
\partial_t
H(t,\th)=\partial^2_{\th}H(t,\th)-\left(\mu+\frac{2c(2c-1)}{\sn^2\th}\right)
(\sn^2\th)^c\dn\th P_{00}(t)\eeq
with
$$H(0,\th)=\dn\th(\sn^2\th)^c.$$
Let us stress that since the moment problem for the associated \SC polynomials
is
{\em determined} the measure $\Psi$ is unique and therefore the solution of the
Kolmogorov equations is unique \cite{re57}.

In order to get it we extend $H(t,\th)$, a priori defined for $\th\in[0,K]$, to
the
interval $\th\in[-K,+K]$ by using the symmetry $\th\leftrightarrow-\th$ of the
equation
and of the boundary value, and further to all real values of $\th$ . This last
step is
possible since (\ref{e6}) and the boundary conditions are periodic with period
$2K$.

We shall use a Laplace transform in the variable t
$$H(t,\th)\longrightarrow\tilde{H}(p,\th)=\int_0^{+\nf}dt\ e^{-pt}H(t,\th)$$
to solve equation (\ref{e6}). We shall first examine what can be said on
general grounds
on $\tilde{H}(p,\th)$.

Firstly since the ${\cal P}_{mn}$ are probabilities we have the bounds
$$\dst 0\leq P_{mn}(t)\leq\frac{(1+c)_n}{(1/2+c)_n}\hspace{1cm}n=0,1\cdots
\hspace{1cm}t\geq 0$$
which imply
$$0\leq H(t,\th)\leq\dn\th(\sn^2\th)^c\fx{2}{1}{1,1+c}{1/2+c}{\sn^2\th}$$
for any real $\th$. From theorem 2.1 of \cite[p.38]{wid46} it follows that
$\tilde{H}(p,\th)$ is analytic in the domain ${\rm Re}\ p>0$ uniformly for real
$\th$.

Secondly the small time behaviour of the transition probabilities is given by
$$\lim_{t\to 0}{\cal P}_{mn}(t)=\de_{mn}$$
from this and theorem 1 of \cite[p.181]{wid46} we conclude to
$$\lim_{p\to+\nf}\tilde{H}(p,\th)=\lim_{t\to 0}H(t,\th)=\dn\th(\sn^2\th)^c.$$

Thirdly $\tilde{H}(p,\th)$ must be periodic in $\th$, with period 2K, and for
$P_{00}(t)$
to exist it is necessary that
\beq\label{e7}
\lim_{\th\to 0}H(t,\th)=0\hspace{1cm}t\geq 0.\eeq

Having stated the most useful properties of $\tilde{H}(p,\th)$ let us now take
the Laplace
transform of equation (\ref{e6}). We get
$$\partial^2_{\th}\tilde{H}(p,\th)-p\tilde{H}(p,\th)=
-H(0,\th)+A(\th)\tilde{P}_{00}(p)=J(p,\th)$$
with
$$A(\th)=\left(\mu+\frac{2c(2c-1)}{\sn^2\th}\right)\dn\th(\sn^2\th)^c.$$
This equation has for general solution even in $\th$
$$\dst\tilde{H}(p,\th)=\frac{e^{\sqrt{p}\th}}{2\sqrt{p}}
\int_0^{\th}d\phi\ e^{-\sqrt{p}\phi}J(p,\phi)+\left(\th\leftrightarrow
-\th\right)
+C(p)\cosh(\sqrt{p}\th)$$
The integral over $\phi$ is convergent at 0 provided that we take $c>1/2$.

The necessary condition (\ref{e7}) implies $C(p)=0$ and the periodicity of
$\tilde{H}(p,\th)$ in the variable $\th$ implies
$$\dst\int_0^{2K}d\phi\ e^{-\sqrt{p}\phi}J(p,\phi)=0\hspace{1cm}c>1/2$$
from which we deduce
$$\tilde{P}_{00}(p)=\dst\frac{\int_0^{2K}d\phi\ e^{-\sqrt{p}\phi}\ H(0,\phi)}
{\int_0^{2K}d\phi\ e^{-\sqrt{p}\phi}\ A(\phi)}.$$
Using the periodicity of A and $\tilde{H}$ brings this ratio to
$$\dst\tilde{P}_{00}(p)=\frac{\int_0^{K}d\phi\cosh{\sqrt{p}(K-\phi)}\
H(0,\phi)}
{\int_0^{K}d\phi\cosh{\sqrt{p}(K-\phi)}\ A(\phi)}$$
and the change of variable $p=-z$ gives eventually the Stieltjes transform
$$\dst\int_0^{+\nf}\frac{d\Psi}{z-s}=-\frac{\int_0^{K}du\ \cos{\sqrt{z}(K-u)}\
H(0,u)}
{\int_0^{K}du\ \cos{\sqrt{z}(K-u)}\ A(u)}$$
so that if we define
$$D(c,\mu;z)=\dst\int_0^{K}du\ \cos{\sqrt{z}(K-u)}\left(2c(2c-1)+\mu\,\sn^2
u\right)
\dn u\,\frac{(\sn^2 u)^{c-1}}{\Ga(2c+1)}$$
we end up with
\beq\label{e8}
\dst\int_0^{+\nf}\frac{d\Psi(s)}{z-s}=-\frac{D(c+1,0;z)}{D(c,\mu;z)}
\hspace{1cm}c>1/2\eeq
in perfect agreement with the result derived in \cite{val94}. In this reference
the
polynomials with rates
$$\la_n=(2n+2c+1)^2,\hspace{2cm}\mu_n=4k^2(n+c)^2+\mu k^2\de_{n0}
\hspace{1cm}0<k^2<1$$
have been shown to follow from the result (\ref{e8}) using the transformation
theory for
Jacobian elliptic functions.

 \end{document}